\documentstyle[12pt,twoside, epsf]{article}
\newcommand{\RR}{\mbox{$I\!\! R$}}	% die reellen Zahlen
 \newtheorem{th}{Theorem}
 \newtheorem{lem}{Lemma}

 \newtheorem{rem}{Remark}

\def\picill#1by#2(#3)
{\vbox to #2
{\hrule width #1 height 0pt depth 0pt
\vfill\epsffile{#3}}}

\let \ttorg \tt \def \tt{\ttorg \obeyspaces}

\begin{document}

\pagestyle{myheadings}

\date{}

\markboth{{\sc Kauffman \& Lambropoulou}}{{\sc Virtual Braids}}

\title{\bf Virtual Braids}

\author{\sc Louis H. Kauffman and Sofia Lambropoulou}

 \maketitle
  
 \thispagestyle{empty}

\section{Introduction}

Just as classical knots and links can be represented by the closures of braids, so can 
 virtual
knots and links be represented by the closures of virtual braids \cite{VKT}. Virtual
 braids have a
group structure that can be described by generators and relations, generalizing the generators and
relations of the classical braid group. This structure of virtual braids is worth study for its own
sake. The virtual braid group is an extension of the classical braid group by the symmetric group.
In \cite{Kamada} a Markov Theorem is proved for virtual braids, giving a set of moves on virtual
braids that generate the same equivalence classes as the virtual link types of their closures. Such
theorems are important for understanding the structure and classification of virtual knots and
links. 
\bigbreak

In the present paper we give a new method for converting virtual knots and 
links to virtual  braids. Indeed, the braiding method given in this paper is quite general and
applies to all the categories in which braiding can be accomplished. This includes the  braiding of 
classical, virtual, flat, welded, unrestricted, and singular knots and links.
We also give reduced presentations for the virtual braid group and for the flat virtual braid group (as
well as for other categories). These reduced presentations are based on the fact that these virtual braid
groups for
$n$ strands are generated by a single braiding element plus the generators of the symmetric group on $n$
letters.
\bigbreak

In a sequel to this paper we shall give a new way to establish the Markov Theorem for
virtual braids via the $L-$move \cite{LR}. The $L-$move has provided a new approach to
Markov Theorems in classical low dimensional topology, and it performs a similar role for virtual
link theory. We shall recover the Markov Theorem of \cite{Kamada} and put it into the context  of
$L$-move theory.  In that sequel paper we shall discuss the same issues for
welded and flat virtual braids.
\bigbreak

This paper consists in four sections. In Section 1 we review the definitions of virtual knot theory, 
and flat virtual knot theory. We recall the interpretations of virtual knot theory in terms of abstract
Gauss codes, and in terms of stabilized embeddings (immersions for flat virtuals) in thickened surfaces.
We emphasize the role of the detour move, and the role of the forbidden moves in the structure of the
theory. A useful feature of this introduction is our description of ribbon neighborhood representations
for virtual links virtual braids. These representations (also called abstract link diagrams) give the
least surface embedding (with boundary) that can  represent a given link diagram. In Section 2 we give a
new braiding algorithm for virtual braids. This algorithm generalizes Alexander's original algorithm,
converting classical links to braid form. The present algorithm is quite general, and gives a uniform
braiding method for many different categories. Section 3 gives the definition of the virtual braid group
via generators and relations. We see the virtual braid group on $n$ strands, $VB_{n},$ as an extension
of the classical braid group $B_{n}$ by the symmetric group $S_{n}.$ The relationship between $S_{n}$
and $B_{n}$ in $VB_{n}$ is intricate. One remarkable property of these subgroups is that
$VB_{n}$ can be generated by $S_{n}$ and a single generator of $B_{n}$ (a single twist, e.g. the braid 
$\sigma_{1}$).  We give a reduced presentation of $VB_{n}$ that incorporates this reduction. In Section
4 we give a similar reduced presentation for the flat virtual braid group $FV_{n}.$ The flat group
$FV_{n}$ is a free product, with extra relations, of two copies of the symmetric group $S_{n}.$ 
 In Section 5 we detail reduced presentations for the welded braid group, the unrestricted virtual braid
group and the flat unrestricted braid group. Section 6 concludes the paper with a topological
interpretation of welded and flat unrestricted braids in terms of tubes imbedded in four-space.
 All of these topics and formulations will be explored further in our subsequent papers.
\bigbreak

For reference to previous work on virtual braids the reader should consult  \cite{CS1, DK, GPV, HR, HRK,
VKT, SVKT, DVK, KADOKAMI, Kamada, KamadaNS, KiSa, KUP, M, Satoh, TURAEV, V1, V2}. For work on welded braids, see
\cite{FRR, Kamada}. For work on  singular braids, see \cite{Baez, Birman, FRZ, KV, V2}.

\section{Virtual Knot Theory}

Virtual knot theory is an extension of classical diagrammatic knot theory. In this extension one 
adds  a {\em virtual crossing} (See Figure 1) that is neither an over-crossing
nor an under-crossing.  A virtual crossing is represented by two crossing arcs with a small circle
placed around the crossing point.  
\bigbreak

Figure 1 illustrates the simplest example of a virtual knot. Note the appearance of a virtual crossing 
in the diagram. There is no way to represent in the plane the Gauss code of this diagram (shown in the
figure)  without entering a virtual crossing that is not registered in the code itself.
Similar remarks apply to links and virtual links, where the codes are collections of sequences, one for
each component of the link. One way to understand the structure of virtual  knots and links is to regard
their diagrams as planar  representatives of possibly non-planar Gauss codes. The virtual crossings are
artifacts of the planar representation.
\bigbreak

\begin{center}
{\tt    \setlength{\unitlength}{0.92pt}
\begin{picture}(148,201)
\thicklines   \put(12,13){$o1+u2+u1+o2+$}
              \put(84,154){$2$}
              \put(17,154){$1$}
              \put(39,189){\vector(1,0){43}}
              \put(38,154){\vector(0,1){36}}
              \put(89,148){\vector(1,0){47}}
              \put(82,34){\vector(-1,0){70}}
              \put(82,190){\vector(0,-1){156}}
              \put(37,78){\vector(0,1){62}}
              \put(136,78){\vector(-1,0){99}}
              \put(136,147){\vector(0,-1){68}}
              \put(13,148){\vector(1,0){59}}
              \put(12,36){\vector(0,1){112}}
              \put(82,78){\circle{22}}
\end{picture}}

{ \bf Figure 1 -- Virtual Knot and Gauss Code} 
\end{center}  

Moves on virtual diagrams generalize the Reidemeister moves for classical knot and link
diagrams.  See Figures 2 and 3.  One can summarize the moves on virtual diagrams as follows:  
The classical crossings interact with one another according to the usual Reidemeister moves (Part A
of Figure 2). The first move of Part A is called {\it planar isotopy move}. Virtual crossings
interact with one another by Reidemeister moves that ignore the  structure of under or over
crossings (Part B of Figure 2). The key move between virtual and classical crossings is shown in 
Part C of Figure 2.  Here a consecutive sequence of two virtual crossings
can be moved across a single classical crossing. The  moves containing virtual
crossings (the moves in Part B and C of Figure 2) shall be called {\it virtual Reidemeister moves}.
All diagrammatic moves of Figure 2 are called {\it augmented Reidemeister moves} and they give rise
to an equivalence relation in the set of virtual knot and link diagrams, called  {\it virtual
equivalence\/} or {\it virtual isotopy\/} or just {\it isotopy\/}.           

\bigbreak

The move in  Part C of Figure 2 is a special case of the more general {\em detour move\/} indicated
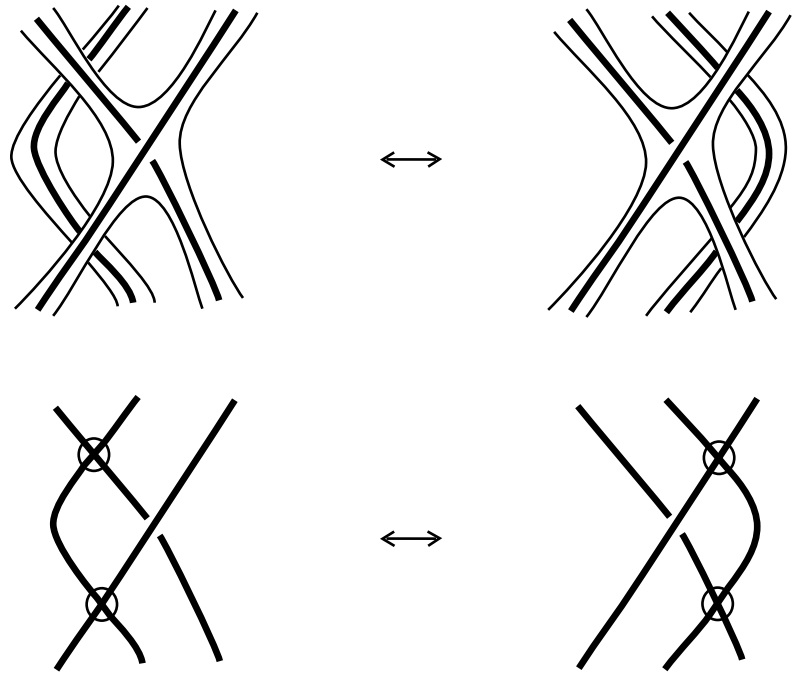
in Figure 3. We will call it the {\it special detour move.} In the detour move, an arc in the diagram
that contains a consecutive sequence of virtual crossings can be excised, and the arc re-drawn,
transversal to the rest of the  diagram (or itself), adding virtual crossings whenever these
intersections occur.  In fact, each of the moves in Parts B and C of Figure~2 can be regarded as
special cases of  the detour move of Figure 3. 

\bigbreak

$$ \picill5inby3.7in(VP2)  $$

\begin{center} 
{ \bf Figure 2 -- Augmented Reidemeister Moves for Virtuals} 
\end{center}  

The equivalence relation generated on virtual diagrams by virtual Reidemeister moves is the same
as the equivalence relation on virtual diagrams generated by the detour move. To see this:
Obviously the moves in B and C of Figure 2 are special cases of the detour move. On the other hand,
by similar arguments as in the classical Reidemeister Theorem, it follows that any detour move can
be achieved by a finite sequence of local steps, each one being a virtual Reidemeister move. Thus
the general detour move is itself the  consequence of the collection of moves in Parts B and C of
Figure 2.
\bigbreak

A succinct description of virtual
equivalence is that it is generated by classical Reidemeister moves {\em and\/} the detour move. 

\bigbreak
 
$$ \picill5inby2.2in(VP3)  $$

\begin{center} 
{ \bf Figure 3 -- The  Detour Move} 
\end{center}

We note that a move analogous to the move in  Part C of Figure 2 but with two real crossings and
one virtual crossing is a {\it forbidden move\/} in virtual knot theory.  There are two types of
forbidden moves: One with an over arc, denoted $F_1,$ and another with an under arc, denoted $F_2.$
 See \cite{VKT} for explanations and interpretations. Variants of the forbidden moves are illustrated in
Figure 4.

\bigbreak
 
$$ \picill5inby1.4in(VP4)  $$

\begin{center} 
{ \bf Figure 4 -- The Forbidden Moves} 
\end{center}

 We know
\cite{VKT,GPV} that classical knot theory embeds faithfully in virtual knot theory. That is, if two 
classical knots are equivalent through moves using virtual crossings, then they are equivalent as
classical knots via standard  Reidemeister moves.

\bigbreak
One can generalize many 
structures in classical knot theory to the virtual domain, and use the virtual knots to test the
limits of classical problems such as  the question whether the Jones polynomial detects knots.  
Counterexamples to this conjecture exist in the virtual domain. It is an
open problem whether some of these counterexamples are isotopic to  classical knots and links.  
\bigbreak

There is a useful topological interpretation for  virtual knot theory in terms of embeddings of 
links in thickened surfaces. Regard each 
virtual crossing as a shorthand for a detour of one of the arcs in the crossing through an 1-handle
that has been attached to the 2-sphere of the original diagram.  
By interpreting each virtual crossing in this way, we
obtain an embedding of a collection of circles into a thickened surface  $S_{g} \times I$ where $g$ 
is the  number of virtual crossings in the original diagram $L$, $S_{g}$ is a compact oriented
surface of genus $g$ and $I$ denotes the unit interval.  We say that two such surface embeddings are
{\em stably equivalent} if one can be obtained from another by isotopy in the thickened surfaces, 
homeomorphisms of the surfaces and the addition or subtraction of empty handles.  Then we have the
following.

\smallbreak
\noindent
{\bf Theorem \cite{VKT,DVK,CS1}.} {\em Two virtual link diagrams are equivalent if and only if
their  correspondent  surface embeddings are stably equivalent.}  

\smallbreak  

See 1.2 below for more discussion of surfaces and virtuals.

\subsection{Flat Virtual Knots and Links}

Every classical knot or link diagram can be regarded as an immersion of cirlces in the plane with extra
structure  at the double points. This extra structure is usually indicated by the over and under crossing
conventions that give instructions for constructing an embedding of the link in three dimensional space
from the diagram.  If we take the diagram without this extra structure, it is the shadow of some link in
three dimensional space, but the weaving of that link is not specified. We call these shadow crossings
{\it flat crossings.} It is well known that if one is
allowed to apply the Reidemeister moves to such a shadow (without regard to the types of crossing since
they are not specified) then the shadow can be reduced to a disjoint union of circles. This reduction is
no longer true in the presence of virtual crossings. 
\bigbreak

More precisely, let a {\em flat virtual diagram} be a diagram with
virtual crossings  and  flat crossings. Two flat virtual diagrams are {\em equivalent} if there is a
sequence of generalized flat Reidemeister moves  taking one to the other. A
{\it generalized flat Reidemeister move} is any move as shown in Figure 2, but with flat crossings in
place of  classical crossings. Note that in studying flat virtuals the rules for changing virtual
crossings among themselves and the rules for changing flat crossings among themselves are identical.
However, detour moves as in Figure 2C are available only for virtual crossings with respect to flat
crossings and not the other way around. The study of flat virtual knots and links was initiated in \cite{VKT}.
The category of flat virtual knots is identical in structure to what are called {\em virtual strings} by V. Turaev in \cite
{TURAEV}. 
\bigbreak

We shall say that a virtual diagram {\em overlies} a flat diagram if the virtual diagram is obtained 
from the flat diagram by choosing a crossing type for each flat crossing in the virtual diagram. To each
virtual diagram $K$ there is an associated flat diagram $F(K)$ that is obtained by forgetting the extra
structure at the classical crossings in $K.$ Note that if $K$ is equivalent to $K'$ as virtual diagrams,
then $F(K)$ is equivalent to $F(K')$ as flat virtual diagrams. Thus, if we can show that $F(K)$ is not
reducible to a disjoint union of circles, then it will follow that $K$ is a non-trivial virtual link.

\bigbreak

$$\vbox{\picill2inby1.6in(VP5)  }$$ 

\begin{center}
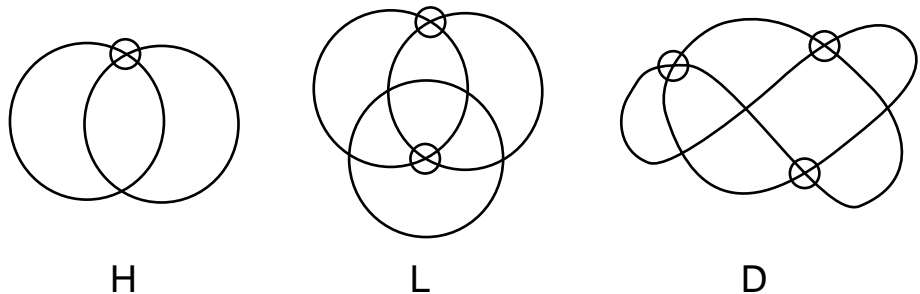
{\bf Figure 5 -- Examples of  Flat Knots and Links} 
\end{center}

Figure 5 illustrates  flat virtual links $H$ and $L$ and a flat virtual knot $D.$ The link $H$ cannot be
undone in the flat  category because it has an odd number of virtual crossings between its two
components and each generalized Reidemeister move preserves the parity of the number of virtual
crossings between components.    The diagram $D$ is shown to be a 
non-trivial flat virtual knot using the filamentation invariant that is introduced in \cite{HR}. The
diagram $L$ is also a non-trivial flat diagram. Note that it comes apart at once if we allow the
forbidden move. 

\bigbreak
The flat virtual diagrams present a strong challenge for the construction of new
invariants.  It is important to understand the structure of flat virtual knots and links. This structure
lies at the heart of the comparison of classical and virtual links. We wish to be able to determine when
a given virtual link is equivalent to a classical link. The reducibility or irreducibility of the
underlying flat diagram is the first obstruction to such an equivalence. See \cite{HR, HRK, TURAEV} for a discussions of combinatorial
invariants of flat virtual knots based on the underlying Gauss code.

\bigbreak
Just as virtual knots and links can be interpreted via stabilized embeddings of curves in thickened
surfaces, flat virtuals can be interpreted as stabilized immersions of curves in surfaces (no thickening
required). See \cite{KADOKAMI} for applications of this point of view.

\subsection{Ribbon Neighborhood Representations}

As we have indicated above, virtual knots and links can be represented as knots and links in thickened 
surfaces. Another way to make this representation is to form a ribbon-neighborhood surface (sometimes
called an {\it abstract link diagram}) \cite{KamadaNS} for a given virtual knot or link, as illustrated in Figure 6. In
this figure we show how a virtual trefoil knot (two classical and one virtual crossing) has the
classical  crossings represented as diagrammatic crossings in disks, which are connected by ribbons,
while the virtual crossings are represented by ribbons that pass over one another without interacting.
The abstract link diagram is shown embedded in three dimensional space, but it is to be regarded as
specified without any particular embedding. Thus it can be represented with the ribbons for the virtual
crossings switched.

\bigbreak
The abstract link diagram is a method for representing a virtual diagram (as an embedding in a thickened
surface) that is distinct from our description in terms of handles given just before 1.1. These two
points of view can be related to one another, and this will be done elsewhere. Here we note that in the
abstract link diagram any closed boundaries can be
filled in with disks or with any convenient orientable surface with boundary. This is in accord with the
representation of virtual knots and links as embeddings in thickened surfaces, taken up to addition and
subtraction of empty handles.
\bigbreak

$$\vbox{\picill1.8inby1.5in(VP6) }$$ 

\begin{center}
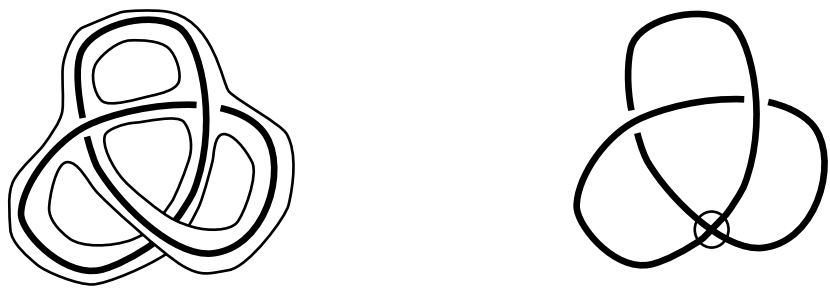
{\bf Figure 6 -- A Ribbon Neighborhood Surface for the Virtual Trefoil} 
\end{center}

In Figure 7 we illustrate the abstract link diagram for one of the special detour moves for virtuals 
(in braided form). Note again how this detour move is accomplished via the freedom of
movement of the virtual crossings represented by non-interacting ribbons.
\bigbreak

$$\vbox{\picill3.5inby3.5in(VP7) }$$ 

\begin{center}
{\bf Figure 7 -- A Special Detour Move} 
\end{center}

In Figure 8 we illustrate a variant of the classical third Reidemeister move in surface form. Note
that we accomplish this move  by adding a disk and then performing an isotopy of the diagram on the
surface.
\bigbreak

$$\vbox{\picill3.5inby3.7in(VP8) }$$ 

\begin{center}
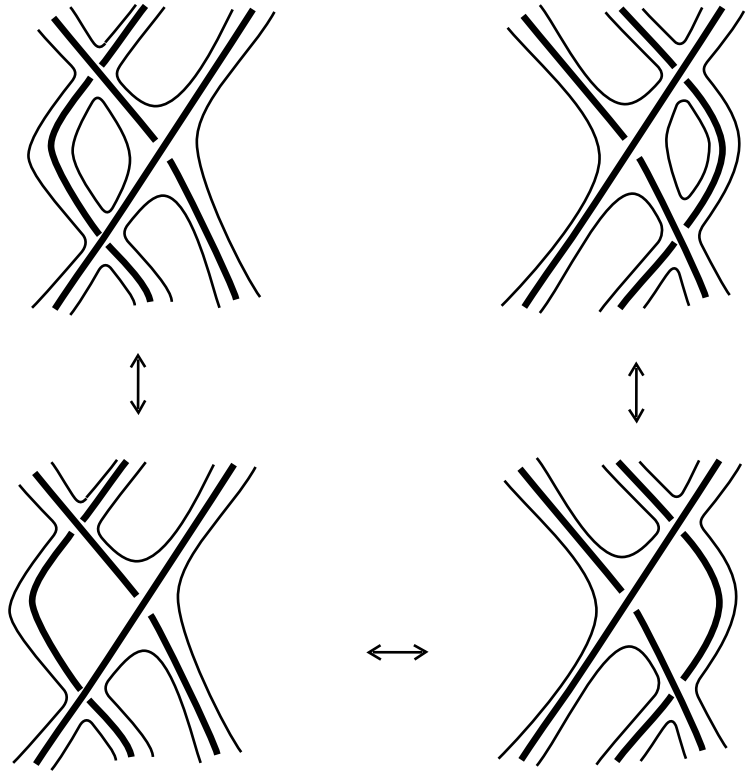
{\bf Figure 8 -- Ribbon Version of Third Reidemeister Move} 
\end{center}

Finally in Figure 9 we illustrate the special detour move for flat virtuals using abstract link diagrams.
Note the stark difference here between the virtual crossing structure and the immersion structure of 
the flat  crossings.
\bigbreak

$$\vbox{\picill3.5inby3.6in(VP9) }$$ 

\begin{center}
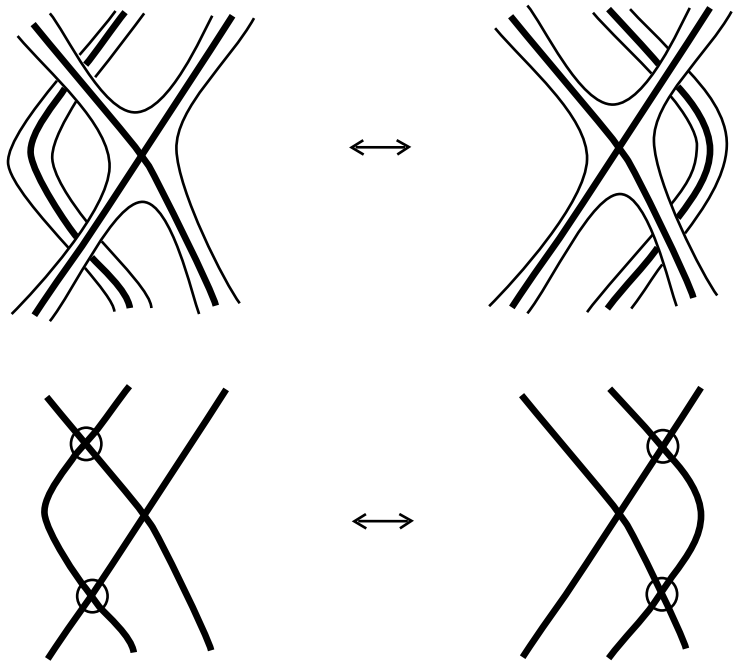
{\bf Figure 9 -- Flat Version of the Detour Move} 
\end{center}

\section{Braiding  Link Diagrams}

A {\it virtual braid} on $n$ strands is a braid on $n$ strands in the classical sense, which may
also contain virtual crossings. The closure of a virtual braid is obtained by joining with simple
arcs the corresponding endpoints of the braid. The set of isotopy classes of virtual braids on $n$
strands forms a group. The virtual braid group structure  will be defined in the next section. In this
section we shall describe a simple and  general method for converting any virtual knot or link diagram 
into the closure of a virtual braid. That is, we shall give a new proof (see \cite{Kamada}) that the
classical Alexander Theorem generalizes to virtuals and, in fact, to all the categories in which braids
are constructed.

\begin{th}{ \ Every 
(oriented) virtual link  can be represented by a virtual braid, whose closure is
isotopic to the original link.
}
\end{th}

\bigbreak

\noindent {\it Proof.} \ The context of our braiding method is based on \cite{LR}. Any  virtual link
diagram can be arranged to be in {\it general position}  with respect to the standard height function on
the plane. This means that it does not contain any horizontal  arcs and it can be seen as a composition
of horizontal stripes, each containing either a local minimum or a local maximum or a crossing (of
classical or virtual type).
\bigbreak

 The idea of the braiding  is on the one
hand to keep the {\it down--arcs} of the diagram that are oriented downwards  and on the other hand to
eliminate the {\it up--arcs}  that go upwards, and produce instead  braid strands. 
First consider  up--arcs that occur between maxima and minima and contain 
no crossings.  Call such an arc in the diagram a {\it free up--arc}. 

\bigbreak
We eliminate a free up--arc  
as follows:  We cut the arc at a point.  We then pull the two arcs the
upper upward and the lower downward, keeping their ends aligned, and so as to obtain a pair of
corresponding braid strands, which create  only  {\it virtual} crossings  with the rest of the diagram.
See Figure 10.  This operation shall be called a {\it braiding move}. 
 The closure of the resulting tangle is  a virtual link diagram, obviously isotopic to the original one.
 Indeed, from the free up--arc we created a stretched loop, which by the detour move is isotopic to the
up--arc.

\bigbreak 

$$\vbox{\picill2inby3.7in(VP10)  }$$ 

\begin{center}
{ \bf Figure 10 --  The  Braiding Move}
\end{center}

Before performing these braiding moves, we prepare the diagram by rotating all the crossings so that 
any arcs that pass through the crossings are directed downwards. There are two
types of rotation: If only one arc in the crossing goes up, then we rotate the crossing by 90 
degrees. If both arcs of the crossing go up, then we rotate it by 180 degrees.  See Figures 11 and 12.
 These rotations may produce new free up--arcs. After adjusting all the crossings, we then braid all the
free up--arcs. The resulting tangle is the desired virtual braid, the closure of which is   isotopic to
the original diagram. 

\bigbreak 

$$\vbox{\picill8inby3.1in(VP11) }$$ 

\begin{center}
{ \bf Figure 11 -- Full Twist} 
\end{center}

\bigbreak 

$$\vbox{\picill8inby2.3in(VP12) }$$ 

\begin{center}
{ \bf Figure 12 -- Half Twist} 
\end{center}

\noindent The braiding algorithm given above will
braid any  virtual diagram and thus it proves the analogue of the Alexander Theorem for virtual
links. \hfill$\Box$

\bigbreak

In Figures 13 and 14 we illustrate an example. In Figure 13 we show a virtual knot and its
preparation for braiding by crossing rotation. In Figure 14 we braid the arcs of the diagaram prepared
in Figure 13.

\bigbreak 

$$\vbox{\picill8inby2.3in(VP13) }$$ 

\begin{center}
{ \bf Figure 13 -- Prepare Example} 
\end{center}

\bigbreak 

$$\vbox{\picill8inby2.7in(VP14) }$$ 

\begin{center}
{ \bf Figure 14 -- Braid Example} 
\end{center}

\begin{rem}{ \rm \  The braiding technique, described in this section, applies equally well to {\it flat 
virtual braids\/} with no change in the description of the procedure. (See Section 4 for definition and
further discussion on the flat braid group.)  For {\it welded\/} and {\it
unrestricted virtual  braids\/} the procedure is the same with welded crossings replacing the role of
virtual crossings. (See Section 5 for definition and further discussion on the corresponding braid
group structures.)  To use this procedure to create a {\it classical braid\/} from a classical knot or
link diagram, the braiding of arcs must be done so that the new braid strands run entirely {\it over\/}
the previously constructed tangle or entirely {\it under\/} it. The same procedure applies also for {\it
singular braids\/}
\cite{Baez, Birman} with singular crossings replacing the role of classical crossings. In this way, this
braiding method works in all the categories in which braids are constructed.  }
\end{rem}

\section{A Reduced Presentation for the Virtual Braid Group}

The set of isotopy classes of virtual braids on $n$ strands forms a group, the {\it virtual braid
group}  denoted $VB_n,$ that was introduced in \cite{VKT}.  The group operation is the usual braid 
multiplication (form $bb'$ by attaching the bottom strand ends of $b$ to the top strand ends of $b'$).  
$VB_n$  is generated by the usual braid generators
$\sigma_{1},\ldots, \sigma_{n-1}$ and by the virtual generators $v_{1},\ldots, v_{n-1},$ 
where each virtual  crossing $v_{i}$ has the form of the braid generator $\sigma_{i}$  with the
crossing replaced by a virtual crossing. See Figure 15 for illustrations. Recall that in  virtual
crossings we do not distinguish between under and over crossing. Thus, $VB_n$  is an extension of
the classical braid group $B_n$ by  the symmetric group $S_n.$  

\bigbreak

 $$\vbox{\picill10cmby1.5in(VP15)  }$$

\begin{center}
{\bf Figure 15 -- The Generators of $VB_n$}
\end{center}
\vspace{3mm}

 Among themselves the braid generators satisfy the usual braiding relations:

$$ \begin{array}{cccl} 
\sigma_{i}\sigma_{i+1}\sigma_{i} & = & \sigma_{i+1}\sigma_{i} \sigma_{i+1}, &   \\
\sigma_i \sigma_j & = & \sigma_j \sigma_i, & \mbox{for} \ j\neq i\pm 1.  \\
\end{array}$$

\noindent Among themselves, the virtual generators are a presentation for the  group $S_{n},$ so
they satisfy the following {\it virtual relations:}  

$$ \begin{array}{cccl} 
{v_{i}}^2 & = & 1, &   \\
v_{i}v_{i+1}v_{i} & = & v_{i+1}v_{i} v_{i+1}, &   \\
v_i v_j & = & v_j v_i, & \mbox{for} \ j\neq i\pm 1.  \\
\end{array}$$

\noindent It is worth noting at this point that the virtual braid group $VB_n$ does not embed in the
classical braid group $B_n,$ since the virtual braid group contains torsion elements and it is
well-known that $B_n$ doesn't. 
 The {\it mixed relations} between
virtual generators and braiding generators are as follows:
 
$$ \begin{array}{cccl} 
\sigma_i v_j & = & v_j \sigma_i, & \mbox{for} \ j\neq i\pm 1,   \\
v_{i} \sigma_{i+1} v_{i} & = & v_{i+1} \sigma_i v_{i+1}. &   \\ 
\end{array}$$

\noindent The second mixed relation will be called the {\it special detour relation} and it is
illustrated in Figure 7. Note that the following  relations are also special detour moves for virtual
braids and they are easy consequences of the above.

$$ \begin{array}{cccl} 
{\sigma_i}^{-1}v_{i+1}v_{i} & = & v_{i+1}v_{i} {\sigma_{i+1}}^{-1}, &   \\
v_{i}v_{i+1} {\sigma_i}^{\pm} & = & {\sigma_{i+1}}^{\pm}v_{i}v_{i+1},&   \\
{\sigma_{i}}^{\pm} v_{i+1}v_{i} & = & v_{i+1}v_{i} {\sigma_{i+1}}^{\pm}.  &  \\
\end{array}$$

\bigbreak

This set of relations taken together define the basic isotopies for virtual braids. Note that each 
relation is a braided version of a local virtual link isotopy.  The special detour relation  is written
equivalently: $\sigma_{i+1} = v_i v_{i+1}
\sigma_i v_{i+1} v_i.$ Notice that this relation is  the braid detour move of the $i$th strand around
the crossing between the $(i+1)$st and the
$(i+2)$nd strand (see illustration in Figure 16) and it provides an inductive way of  expressing
all  braiding generators in terms of the first braiding generator $\sigma_{1}$  and the virtual 
generators $v_1,\ldots,v_{n-1}.$

In this section we give a reduced presentation  for $VB_n$ with generators 
$$\{\sigma_1, v_1,\ldots,v_{n-1}\},$$ 
where we assume the   defining relations:

$$\sigma_{i+1} := (v_i\ldots v_2v_1)\,(v_{i+1}\ldots v_3v_2)\,\sigma_1\, (v_2v_3\ldots
v_{i+1})\,(v_1v_2\ldots v_i)  \  \  \ (*)$$

\bigbreak

\noindent for \  $i=1,\ldots, n-2.$ In terms of braid diagrams, this relation is the braid detour
moves of the strands $1, 2, \ldots, i$  around the crossing $\sigma_{i+1}.$ See right hand illustration
in Figure 16.  

\bigbreak

 $$\vbox{\picill15cmby2in(VP16)  }$$

\begin{center}
{\bf Figure 16 -- Detouring the Crossing $\sigma_{i+1}$}
\end{center}
\vspace{3mm}

\begin{rem}{ \rm \  By the detour move, in the same way that a crossing can be detoured to the first
position of the braid, similarly any box in the braid can be detoured to the first position (in
fact, to any position), see Figure 17. There it may undergo some  changes and then it can be detoured
back to its original position in the braid.  In particlular, a relation in $VB_n$ occuring in a box in
the interior of a braid may be redundant.  In order to omit it we detour the box to the first
position, there we apply a specific relations (see statement of Theorem 2) and then we detour the result
back, thus obtaining the other side of the relation that we wanted to omit.
 }
\end{rem}

\bigbreak

 $$\vbox{\picill11cmby2.7in(VP17)  }$$

\begin{center}
{\bf Figure 17 -- Detouring a Box}
\end{center}
\vspace{3mm}

\begin{th}{ \ The virtual braid group $VB_{n}$ has the following reduced presentation.  

\[VB_{n} = \left< \begin{array}{ll}  
\begin{array}{l} 
\sigma_1, v_1, \ldots ,v_{n-1}  \\
\end{array} &
\left|
\begin{array}{l} 

v_i v_{i+1} v_i = v_{i+1} v_i v_{i+1}    \\

v_i v_j=v_j v_i, \ \  \ j\neq i\pm 1   \\ 

{v_i}^2  = 1, \ \ \ 1 \leq i \leq n-1  \\

\sigma_1 v_j=v_j \sigma_1, \ \ \ j>2   \\ 

 (v_1 \sigma_1 v_1)\,(v_2 \sigma_1 v_2)\, (v_1 \sigma_1 v_1) =
  (v_2 \sigma_1 v_2)\, (v_1 \sigma_1 v_1) \, (v_2 \sigma_1 v_2)   \\ 

\sigma_1\, (v_2 v_3 v_1 v_2 \sigma_1 v_2 v_1 v_3 v_2) =
    (v_2 v_3 v_1 v_2 \sigma_1 v_2 v_1 v_3 v_2) \, \sigma_1   \\ 
\end{array} \right.  \end{array} \right>.  \]
}
\end{th}

In Figure 18  we illustrate the last two  reduced relations. 

\bigbreak

 $$\vbox{\picill12cmby2.5in(VP18)  }$$

\begin{center}
{\bf Figure 18 -- The Main Reduced Relations}
\end{center}
\vspace{3mm}

\noindent Note that the relation: 
$$(v_1 \sigma_1 v_1)\,(v_2 \sigma_1 v_2)\, (v_1 \sigma_1 v_1) =
  (v_2 \sigma_1 v_2)\, (v_1 \sigma_1 v_1) \, (v_2 \sigma_1 v_2)$$ 

\noindent is equivalent to the relation:
$$\sigma_1 \,(v_1 v_2 \sigma_1 v_2 v_1)\, \sigma_1 =
    (v_1 v_2 \sigma_1 v_2 v_1) \, \sigma_1\, (v_1 v_2 \sigma_1 v_2 v_1),$$

 \noindent    which reflects the braid relation: 
$\sigma_1 \sigma_2 \sigma_1 = \sigma_2 \sigma_1 \sigma_2.$ Also, that the relation: 
$$\sigma_1\, (v_2 v_3 v_1 v_2 \sigma_1 v_2 v_1 v_3 v_2) =
    (v_2 v_3 v_1 v_2 \sigma_1 v_2 v_1 v_3 v_2) \, \sigma_1$$

\noindent reflects the braid relation: $\sigma_1 \sigma_3 = \sigma_3 \sigma_1.$  Therefore:

\bigbreak
\noindent \begin{itemize}
\item From the commuting relations: $\sigma_i \sigma_j = \sigma_j \sigma_i$
    we only need to keep: $\sigma_1 \sigma_3 = \sigma_3 \sigma_1.$
\item From the  relations: $\sigma_i  \sigma_{i+1} \sigma_i = \sigma_{i+1} \sigma_i \sigma_{i+1}$
    we only need to keep: $\sigma_1 \sigma_2 \sigma_1  = \sigma_2 \sigma_1 \sigma_2.$
\end{itemize}
\bigbreak

\noindent {\it Proof of Theorem 2.} \ We note first that the special detour relations  
$v_{i} \sigma_{i+1} v_{i}  =  v_{i+1} \sigma_i v_{i+1}$ can be omitted  in the reduced presentation,
since they were used in the defining new relations (*) for the $\sigma_j$'s. The proof of the reduced
presentation is then a consequence of the  three lemmas below.  The proofs of these lemmas are
based on the simple geometric idea described in Remark 2. In
the proofs  we underline  the expressions  that are used each time in the next step. 

\begin{lem}{ \ The mixed relations $\sigma_i v_j  =  v_j \sigma_i$ \  for $i>1$ and $j\neq i\pm
1$    follow from the defining relations $(*),$ the virtual relations and the reduced relations
$\sigma_1 v_j  =  v_j \sigma_1,$ \  for  $j>2.$
 }
\end{lem}

\noindent {\em Proof.} \ Substituting we have: \ $\sigma_i v_j  \stackrel{(*)}{=} (v_{i-1}\ldots
v_2v_1)\,(v_i\ldots v_3v_2)\,\sigma_1\, (v_2v_3\ldots v_i)\,(v_1v_2\ldots v_{i-1})\, v_j.$  Since 
 $j\neq i\pm 1$ either $j\geq i+2$ or $j\leq i-2.$   If $j\geq i+2,$ then in the above
expression $v_j$ clearly commutes  with all generators, thus $\sigma_i v_j  =  v_j \sigma_i.$ 
 If $j\leq i-2$ we have:

\vspace{.1in}
\noindent $\sigma_i v_j  =  (v_{i-1}\ldots
v_2v_1)\,(v_i\ldots v_3v_2)\,\sigma_1\, (v_2v_3\ldots v_i)\,(v_1v_2\ldots v_{i-1})\, 
\underline{v_j} $

\vspace{.1in}
$  =  (v_{i-1}\ldots v_2v_1)\,(v_i\ldots
v_3v_2)\,\sigma_1\, (v_2v_3\ldots v_i)\,(v_1v_2 \ldots v_{j-1}\underline{v_j v_{j+1} v_j}  v_{j+2}
\ldots v_{i-1}) $

\vspace{.1in}
$  =  (v_{i-1}\ldots v_2v_1)\,(v_i\ldots v_3v_2)\,\sigma_1\, (v_2v_3\ldots v_{j+1} v_{j+2} \ldots
v_i)\,(v_1v_2 \ldots v_{j-1} \underline{v_{j+1}} v_j v_{j+1} v_{j+2} \ldots v_{i-1}) $

\vspace{.1in}
$  =  (v_{i-1}\ldots v_2v_1)\,(v_i\ldots v_3v_2)\,\sigma_1\, (v_2v_3\ldots  v_j \underline{
v_{j+1} v_{j+2} v_{j+1}}  v_{j+3} \ldots v_i) \,(v_1v_2 \ldots v_{i-1}) $

\vspace{.1in}
$  =  (v_{i-1}\ldots v_2v_1)\,(v_i\ldots v_{j+2} v_{j+1} v_j \ldots 
v_3v_2)\,\sigma_1\, (v_2v_3\ldots  v_j\underline{ v_{j+2}}
v_{j+1} v_{j+2} v_{j+3} \ldots v_i) \,(v_1v_2 \ldots v_{i-1}) $

\vspace{.1in}
$  =  (v_{i-1}\ldots v_2v_1)\,(v_i\ldots v_{j+3} \underline{ v_{j+2}  v_{j+1} v_{j+2}} v_j\ldots
v_3v_2) \,\sigma_1\, (v_2v_3\ldots  v_i) \,(v_1v_2 \ldots v_{i-1}) $

\vspace{.1in}
$  =  (v_{i-1}\ldots  v_{j+1} v_j v_{j-1} \ldots v_2v_1)\,(v_i\ldots  v_{j+3} \underline{ v_{j+1}} 
v_{j+2} v_{j+1} v_j\ldots v_3v_2)
\,\sigma_1\, (v_2v_3\ldots  v_i) \,(v_1v_2 \ldots v_{i-1}) $

\vspace{.1in}
$  =  (v_{i-1}\ldots v_{j+2} \underline{v_{j+1} v_j v_{j+1} } v_{j-1} \ldots
v_2v_1)\,  (v_i\ldots v_3v_2)
\,\sigma_1\, (v_2v_3\ldots  v_i) \,(v_1v_2 \ldots v_{i-1}) $

\vspace{.1in}
$  =  (v_{i-1}\ldots v_{j+2} \underline{v_j }  v_{j+1}  v_j  v_{j-1} \ldots v_2v_1)\,  (v_i\ldots
v_3v_2) \,\sigma_1\, (v_2v_3\ldots  v_i) \,(v_1v_2 \ldots v_{i-1}) $

\vspace{.1in}
$  =  v_j\, (v_{i-1}\ldots  v_2v_1)\,  (v_i\ldots v_3v_2) \,\sigma_1\,
(v_2v_3\ldots  v_i) \,(v_1v_2 \ldots v_{i-1}) $

\vspace{.1in}
$  \stackrel{(*)}{=}  v_j \sigma_i. $ \hfill $\Box$

\bigbreak

\noindent In the proofs of Lemmas 2 and  3 below we use repeatedly the following virtual braid
relations, which are   easy consequences of the virtual  relations.

$$ v_iv_{i-1}\ldots v_{j+1}{v_j} v_{j+1}\ldots v_{i-1} v_i =
v_jv_{j+1}\ldots v_{i-1}{v_i} v_{i-1}\ldots v_{j+1}{v_j} \ \ \ (\dagger)$$ 

$$(v_4 v_3 v_2 v_1) \ldots (v_{i+2} v_{i+1} v_i v_{i-1}) =  (v_4 \ldots v_{i+2})\,
(v_3 \ldots v_{i+1})\, (v_2 \ldots v_i) \, (v_1 \ldots v_{i-1}). \ \ \ (\ddagger)$$

\begin{lem}{ \ The braid relations \ 
$\sigma_i  \sigma_{i+1} \sigma_i = \sigma_{i+1} \sigma_i
\sigma_{i+1}$
\  for   $i>1$ follow from the defining relations $(*),$ the virtual relations, the reduced
relations $\sigma_1 v_j  =  v_j \sigma_1$ of Lemma~1 
 and  the reduced relation:
$$\sigma_1 \,(v_1 v_2 \sigma_1 v_2 v_1)\, \sigma_1 =
    (v_1 v_2 \sigma_1 v_2 v_1) \, \sigma_1\, (v_1 v_2 \sigma_1 v_2 v_1).$$
 }
\end{lem}

\noindent {\em Proof.} \ Indeed, on the one hand we have: 

\vspace{.1in}
\noindent $\sigma_i  \sigma_{i+1} \sigma_i \stackrel{(*)}{=}   $

\vspace{.1in}
 $= [(v_{i-1}\ldots v_1)\,(v_i\ldots v_2)\,\sigma_1\, (v_2\ldots
v_i)\,(\underline{v_1\ldots v_{i-1})]\,[(v_i v_{i-1}\ldots v_1})\,(v_{i+1} \ldots v_2)\, \sigma_1
\cdot$

\vspace{.1in}
$ \ \ \ \  (v_2\ldots v_{i+1})\,(\underline{v_1 \ldots v_{i-1} v_i)]\,
[(v_{i-1}\ldots v_1}) \,(v_i\ldots v_2)\,\sigma_1\,
(v_2 \ldots v_i)\,(v_1 \ldots v_{i-1})]$

\vspace{.1in}
 $ \stackrel{(\dagger)}{=} (v_{i-1}\ldots v_1)\,(v_i\ldots v_2)\,\sigma_1\, 
(\underline{v_2 \ldots v_i)\, (v_i \ldots v_2} v_1 v_2 \ldots v_i)
 \,(v_{i+1} \ldots v_2)\, \sigma_1 \cdot$

\vspace{.1in}
$ \ \ \   (v_2\ldots  v_{i+1})\, (v_i \ldots v_2 v_1 \underline{v_2 \ldots v_i)\, ( v_i \ldots 
v_2}) \, \sigma_1\,  (v_2\ldots v_i)\,(v_1\ldots v_{i-1}) $

\vspace{.1in}
 $= (v_{i-1}\ldots v_1)\,(v_i\ldots v_2)\,\sigma_1\, 
(v_1 \underline{ v_2\ldots v_i)\, (v_{i+1} v_i\ldots v_2})\,\sigma_1  (\underline{v_2 \ldots
v_i v_{i+1})\, (v_i \ldots v_2} v_1) \, \sigma_1 \cdot$

\vspace{.1in}
$ \ \ \  (v_2\ldots v_i)\,(v_1\ldots v_{i-1}) $

\vspace{.1in}
 $\stackrel{(\dagger)}{=}  (v_{i-1}\ldots v_1)\,(v_i\ldots v_2)\,\sigma_1 
v_1 \,  (\underline{v_{i+1} \ldots v_3} v_2 v_3\ldots v_{i+1})\,\underline{\sigma_1} (v_{i+1} \ldots
v_3v_2 v_3\ldots v_{i+1})\, \underline{v_1\sigma_1} \cdot$

\vspace{.1in}
$ \ \ \   (v_2\ldots v_i)\,(v_1\ldots v_{i-1}) $

\vspace{.1in}
 $= (v_{i-1}\ldots v_1)\,(v_i\ldots v_2)\, (v_{i+1} \ldots v_3)\, \sigma_1 v_1 v_2\sigma_1\,
(\underline{ v_3\ldots  v_{i+1})\, (v_{i+1} \ldots v_3}) \, v_2 v_1 \sigma_1\,  (v_3 \ldots 
v_{i+1}) \cdot  $

\vspace{.1in}
$ \ \ \ (v_2\ldots v_i)\,(v_1\ldots v_{i-1}) $

\vspace{.1in}
 $= (v_{i-1}\ldots  v_1)\,(v_i\ldots  v_2)\, (v_{i+1} \ldots  v_3)\,
\underline{\sigma_1 v_1 v_2 \sigma_1 v_2 v_1 \sigma_1 } \, 
 (v_3 \ldots  v_{i+1}) \, (v_2\ldots v_i) \, (v_1\ldots v_{i-1})$

\vspace{.1in}
 $= (v_{i-1}\ldots v_1)\,(v_i\ldots  v_2)\, (v_{i+1} \ldots  v_3)\,
\underline{ v_1} v_2 \sigma_1 v_2 v_1 \sigma_1 v_1 v_2 \sigma_1 v_2 \underline{ v_1} \, 
 (v_3 \ldots  v_{i+1}) \, (v_2\ldots v_i) \cdot$

\vspace{.1in}
 $ \ \ \ (v_1\ldots v_{i-1})$

\vspace{.1in}
 $= (\underline{v_{i-1}\ldots  v_1)\,(v_i\ldots  v_1)\, (v_{i+1} \ldots  v_2})\,
\sigma_1 v_2 v_1 \sigma_1 v_1 v_2 \sigma_1 \, (\underline{ v_2  \ldots  v_{i+1})  (v_1 \ldots v_i)
\,(v_1 \ldots v_{i-1}})$

\vspace{.1in}
$ = A.$

\bigbreak

\noindent On the other hand with similar manipulations we obtain: 

\vspace{.1in}
\noindent $ \sigma_{i+1} \sigma_i  \sigma_{i+1} =  \cdots $

%\vspace{.1in}
%\noindent $ \sigma_{i+1} \sigma_i  \sigma_{i+1} \stackrel{(*)}{=}   $

%\vspace{.1in}
% $= [(v_i\ldots v_1)\, (v_{i+1}\ldots v_2)\,\sigma_1\, (v_2 \ldots v_{i+1})\,
%(\underline{ v_1\ldots v_{i-1} v_i)\, ( v_{i-1}\ldots v_1)}\, ( v_i\ldots v_2) \, \sigma_1
%\cdot$

%\vspace{.1in}
%$ \ \ \  (v_2\ldots v_i)\, (\underline{v_1\ldots v_{i-1})\,(v_iv_{i-1} \ldots v_1}) \,
%(v_{i+1}\ldots v_2)\, \sigma_1\, (v_2 \ldots v_{i+1})\, (v_1 \ldots v_i)$

%\vspace{.1in}
% $= (v_i\ldots v_1)\, (v_{i+1}\ldots v_2)\,\sigma_1\, (\underline{ v_2 \ldots v_i v_{i+1})\,
%( v_i \ldots v_2} v_1 \underline{ v_2 \ldots v_i)\, ( v_i\ldots v_2}) \, \sigma_1 \cdot$

%\vspace{.1in}
%$ \ \ \  (\underline{ v_2\ldots v_i)\, (v_i \ldots v_2} v_1 \underline{ v_2 \ldots v_i)
%(v_{i+1} v_i \ldots v_2})\, \sigma_1\, (v_2 \ldots v_{i+1})\, (v_1\ldots v_i)$

%\vspace{.1in}
% $= (v_i\ldots v_1)\, (v_{i+1}\ldots v_2)\,\sigma_1 \, 
%(\underline{ v_{i+1}\ldots v_3} v_2 v_3 \ldots v_{i+1})\, v_1 \sigma_1 v_1 \, (v_{i+1}\ldots 
%v_3 v_2 v_3 \ldots v_{i+1})\, \underline{\sigma_1} \cdot$

%\vspace{.1in}
%$ \ \ \  (v_2  \ldots v_{i+1})\, (v_1\ldots v_i)$

%\vspace{.1in}
% $= (v_i\ldots v_1)\, (v_{i+1}\ldots v_2)\,  (v_{i+1}\ldots v_3)\,\sigma_1\, 
%(v_2 v_3 \ldots v_{i+1})\, \underline{v_1 \sigma_1 v_1} \, (v_{i+1}\ldots v_3v_2)\, 
%\sigma_1 \, (v_3 ldots v_{i+1})\cdot$

%\vspace{.1in}
%$ \ \ \   (v_2 \ldots v_{i+1}) \,(v_1\ldots v_i)$

\vspace{.1in}
 $= (\underline{ v_i\ldots v_1)\, (v_{i+1}\ldots v_2)\,  (v_{i+1}\ldots v_3})\,\sigma_1 v_2 v_1
\sigma_1 v_1 v_2
\sigma_1 \, (\underline{ v_3 \ldots v_{i+1})\, (v_2 \ldots v_{i+1}) \,(v_1 \ldots v_i})$

\vspace{.1in}
$ = B.$

\bigbreak 

\noindent But: \, $(v_i\ldots v_1)\, (v_{i+1}v_iv_{i-1} \ldots v_2)\, (\underline{v_{i+1}} v_i
\ldots v_3)$

\vspace{.1in}
$ = (v_i \ldots v_1)\, \underline{v_{i+1}v_i v_{i+1}}\, (v_{i-1} \ldots v_2)\, ( v_i \ldots v_3)$

\vspace{.1in}
$=(v_iv_{i-1} v_{i-2} \ldots v_1)\, \underline{v_i} v_{i+1}v_i\, (v_{i-1} \ldots v_2)\, ( v_i \ldots
v_3)$

\vspace{.1in}
$ =\underline{ v_iv_{i-1} v_i} (v_{i-2} \ldots v_1)\, (v_{i+1} \ldots v_2)\, ( v_i \ldots v_3)$

\vspace{.1in}
$ = v_{i-1}v_i v_{i-1}\, ( v_{i-2} \ldots v_1)\, (\underline{v_{i+1}}v_i \ldots v_2)\, 
( v_i \ldots v_3)$

%\vspace{.1in}
%$ = v_{i-1}v_iv_{i+1}\, (v_{i-1}  \ldots v_1)\, (v_i v_{i-1} v_{i-2} \ldots v_2)\, 
%(\underline{v_i} v_{i-1} \ldots v_3)$

%\vspace{.1in}
%$ = v_{i-1}v_iv_{i+1}\, (v_{i-1} v_{i-2} v_{i-3} \ldots v_1)\, \underline{ v_i v_{i-1} v_i}
%\, ( v_{i-2} \ldots v_2)\, (v_{i-1}\ldots v_3)$

%\vspace{.1in}
%$ = (v_{i-1}v_iv_{i+1})\, \underline{v_{i-1} v_{i-2} v_{i-1}} \, 
%(v_{i-3} \ldots v_1)\, v_i v_{i-1}( v_{i-2} \ldots v_2)\, (v_{i-1}\ldots v_3)$

%\vspace{.1in}
%$ = (v_{i-1}v_iv_{i+1})\,  v_{i-2} v_{i-1} v_{i-2} \, ( v_{i-3} \ldots v_1)\, 
%\underline{ v_i}\, (v_{i-1}\ldots v_2)\, (v_{i-1}\ldots v_3)$

%\vspace{.1in}
%$ = (v_{i-1}v_iv_{i+1})\, ( v_{i-2} v_{i-1} v_i) \, (v_{i-2}  \ldots v_1)\, 
% (v_{i-1}\ldots v_2)\, (v_{i-1}\ldots v_3) = \cdots $

\vspace{.1in}
$ = (v_{i-1}v_iv_{i+1})\, (v_{i-1}  \ldots v_1)\, (v_i \ldots v_2)\, (v_i \ldots v_3) $

$\vdots$

$ = (\underline{v_{i-1}} v_i v_{i+1})\, ( \underline{v_{i-2}} v_{i-1} v_i)  \ldots (
\underline{v_1} v_2 v_3)
\, v_1 v_2 $

\vspace{.1in}
$ =  (v_{i-1} v_{i-2} \ldots  v_1) \, (v_i v_{i-1} \ldots  v_2) \, ( v_{i+1} v_i \ldots 
v_3) \,\underline{v_1} v_2   $

\vspace{.1in}
$ =  (v_{i-1} \ldots  v_1) \, (v_i  \ldots  v_1) \, ( v_{i+1}  \ldots  v_2).    $

\bigbreak

\noindent Hence, and by the symmetry of the underlined expressions  in $A$  and $B,$ we
showed that  $B=A.$ 
   \hfill $\Box$

\begin{lem}{ \ The braid relations \ 
$\sigma_i  \sigma_j = \sigma_j \sigma_i$
\  for $i>1, \ i<j$ and $j\neq i + 1$ follow from the defining relations $(*),$ the virtual
relations, the reduced relations $\sigma_1 v_j  =  v_j \sigma_1$ of Lemma~1 
 and  the reduced relation:
$$\sigma_1\, (v_2 v_3 v_1 v_2 \sigma_1 v_2 v_1 v_3 v_2) =
    (v_2 v_3 v_1 v_2 \sigma_1 v_2 v_1 v_3 v_2) \, \sigma_1.$$
 }
\end{lem}

\noindent {\em Proof.} \ Indeed, on the one hand we have: 

\vspace{.1in}
\noindent $\sigma_i \sigma_j \stackrel{(*)}{=}   [(v_{i-1}\ldots v_1)\,(v_i\ldots
v_2)\,\sigma_1\, (v_2 \ldots v_i)\,(v_1 \ldots v_{i-1})] \, [(\underline{v_{j-1}\ldots v_{i+2}
v_{i+1}} v_i  \ldots v_1)\,(v_j\ldots v_2)\,\sigma_1 \cdot$

\vspace{.1in}
$ \ \ \ (v_2 \ldots v_j)\,(v_1 \ldots v_{j-1})]$

\vspace{.1in}
$= (v_{j-1}\ldots v_{i+2})\, (v_{i-1}\ldots v_1)\,(v_i\ldots v_2)\,\sigma_1\, (v_2 \ldots v_i)\,
 v_{i+1} \, (\underline{ v_1 \ldots v_{i-1}) \,(v_i \ldots v_1})
\, (v_j\ldots v_2)\,\sigma_1 \cdot$

\vspace{.1in}
$ \ \ \   (v_2\ldots v_j)\,(v_1\ldots v_{j-1})$

\vspace{.1in}
$\stackrel{(\dagger)}{=}   (v_{j-1}\ldots v_{i+2})\,(v_{i-1}\ldots v_1)\,(v_i\ldots
v_2)\,\sigma_1\, 
(\underline{ v_2 \ldots v_i v_{i+1})\, (v_i \ldots v_2} v_1 v_2 \ldots v_i) \,(v_j\ldots
v_2)\,\sigma_1 \cdot$

\vspace{.1in}
$ \ \ \    (v_2\ldots v_j)\,(v_1\ldots v_{j-1})$

\vspace{.1in}
$\stackrel{(\dagger)}{=}  (v_{j-1}\ldots v_{i+2})\,(v_{i-1}\ldots v_1)\,(v_i\ldots v_2)\,\sigma_1\, 
(\underline{ v_{i+1} \ldots v_3} v_2  v_3 \ldots v_{i+1})\, (v_1 \ldots v_i)  \cdot$

\vspace{.1in}
$ \ \ \   (\underline{v_j\ldots v_{i+3} v_{i+2}}v_{i+1} \ldots
v_2)\,\sigma_1 \, (v_2\ldots v_j)\,(v_1\ldots v_{j-1})$

\vspace{.1in}
$= (v_{j-1}\ldots v_{i+2})\,(v_j\ldots v_{i+3})\, (v_{i-1}\ldots v_1)\,(v_i\ldots v_2)\,(v_{i+1}
\ldots v_3)\,  \sigma_1\, ( v_2 \ldots v_{i+1}) \, v_{i+2} \cdot$

\vspace{.1in}
$ \ \ \  (v_1 \underline{ v_2 \ldots v_i) \, ( v_{i+1} v_i \ldots v_2})\,\sigma_1\,  (v_2\ldots
v_j)\,(v_1\ldots v_{j-1})$

\vspace{.1in}
$\stackrel{(\dagger)}{=}  (v_{j-1}\ldots v_{i+2})\,(v_j\ldots v_{i+3})\, (v_{i-1}\ldots
v_1)\,(v_i\ldots v_2)\,(v_{i+1}
\ldots v_3)\,  \sigma_1\, ( v_2 \ldots  v_{i+2} ) \cdot$

\vspace{.1in}
$ \ \ \  v_1\, (\underline{ v_{i+1} \ldots v_3} v_2  v_3 \ldots v_{i+1} )\, \underline{\sigma_1}\, 
(v_2\ldots v_j)\,(v_1\ldots v_{j-1})$

\vspace{.1in}
$= (v_{j-1}\ldots v_{i+2})\,(v_j\ldots v_{i+3})\, (v_{i-1}\ldots v_1)\,(v_i\ldots v_2)\,(v_{i+1}
\ldots v_3)\,  \sigma_1 \cdot$

\vspace{.1in}
$ \ \ \  ( v_2 \underline{v_3 \ldots  v_{i+2} )\, ( v_{i+1} \ldots v_3}) v_1 v_2 
\sigma_1\, ( v_3 \ldots v_{i+1})\,  (v_2\ldots v_j)\,(v_1\ldots v_{j-1})$

\vspace{.1in}
$\stackrel{(\dagger)}{=}  (v_{j-1}\ldots v_{i+2})\,(v_j\ldots v_{i+3})\, (v_{i-1}\ldots
v_1)\,(v_i\ldots v_2)\,(v_{i+1}
\ldots v_3)\,  \sigma_1  \cdot$

\vspace{.1in}
$ \ \ \ v_2 \, (\underline{ v_{i+2} \ldots v_4 }v_3 \underline{v_4 \ldots v_{i+2}} ) \,  v_1 v_2 
\sigma_1\,  ( v_3 \ldots v_{i+1})\,  (v_2\ldots v_j)\,(v_1\ldots v_{j-1})$

\vspace{.1in}
$= (v_{j-1}\ldots v_{i+2})\,(v_j\ldots v_{i+3})\, (v_{i-1}\ldots v_1)\,(v_i\ldots v_2)\,(v_{i+1}
\ldots v_3)\, ( v_{i+2} \ldots v_4)\cdot$

\vspace{.1in}
$ \ \ \ \sigma_1 v_2 v_3 v_1 v_2 \sigma_1  \, \underline{ (v_4 \ldots v_{i+2}) \,
( v_3 \ldots v_{i+1})\,  (v_2\ldots v_j)\,(v_1\ldots v_{j-1})}.$

\bigbreak

\noindent But:

\vspace{.1in}
$(v_4 \ldots v_{i+2}) \, ( v_3 \ldots v_{i+1})\,  (v_2 \ldots v_i v_{i+1} \ldots
v_j)\,(v_1 \ldots  v_{i-1} v_i \ldots v_{j-1}) $

\vspace{.1in}
$ \stackrel{(\ddagger)}{=} (v_4 v_3 v_2 v_1) \ldots  (v_{i+2}
\underline{v_{i+1} v_i} v_{i-1})\,  (\underline{ v_{i+1}} v_{i+2} \ldots v_j)\,(v_i \ldots v_{j-1})$

\vspace{.1in}
$ = (v_4 v_3 v_2 v_1) \ldots  (v_{i+1} v_i v_{i-1} v_{i-2})\,\underline{v_i} \,
(v_{i+2} 
 v_{i+1} v_i v_{i-1})\,  ( v_{i+2} \ldots v_j)\,(v_i  \ldots
v_{j-1}) = \cdots$

\vspace{.1in}
$ = v_2\, (v_4 v_3 v_2 v_1) \ldots  (v_{i+1} v_i v_{i-1} v_{i-2})\,
\, (v_{i+2} 
 v_{i+1} \underline{v_i v_{i-1}})\,  ( v_{i+2} \ldots v_j)\,(\underline{v_i } v_{i+1} \ldots
v_{j-1})$

\vspace{.1in}
$ = v_2\, (v_4 v_3 v_2 v_1) \ldots  (v_{i+1} v_i v_{i-1} v_{i-2})\,
\,\underline{v_{i-1} }\, (v_{i+2}  v_{i+1} v_i v_{i-1})\,  ( v_{i+2} \ldots v_j)\,( v_{i+1} \ldots
v_{j-1}) = \cdots$

\vspace{.1in}
$ = v_2 v_1\, (v_4 v_3 v_2 v_1) \ldots  (v_{i+1} v_i v_{i-1} v_{i-2})\,
\, (\underline{ v_{i+2}  v_{i+1}} v_i v_{i-1})\,  ( \underline{v_{i+2}} v_{i+3} \ldots v_j)\,(
v_{i+1} \ldots v_{j-1}) $

\vspace{.1in}
$ = v_2 v_1\, (v_4 v_3 v_2 v_1) \ldots  (v_{i+1} v_i v_{i-1} v_{i-2})\,
\underline{ v_{i+1}} \, (v_{i+2}  v_{i+1} v_i v_{i-1})\,  ( v_{i+3} \ldots v_j)\,(
v_{i+1} \ldots v_{j-1})  = \cdots$

\vspace{.1in}
$ = v_2 v_1 v_3 \, (v_4 v_3 v_2 v_1) \ldots  (v_{i+1} v_i  v_{i-1} v_{i-2})\,
 (v_{i+2} \underline{ v_{i+1} v_i} v_{i-1})\,  ( v_{i+3} \ldots v_j)\,
(\underline{v_{i+1}} v_{i+2} \ldots v_{j-1}) = \cdots$

\vspace{.1in}
$ = v_2 v_1 v_3 v_2\, (v_4 v_3 v_2 v_1) \ldots 
 (v_{i+2}  v_{i+1} v_i v_{i-1})\, ( v_{i+3} \ldots v_j)\,
(v_{i+2} \ldots v_{j-1}).  $

\bigbreak

\noindent Thus: 

\vspace{.1in}
\noindent \ $\sigma_i \sigma_j = (v_{j-1}\ldots v_{i+2})\,(v_j\ldots v_{i+3})\, (v_{i-1}\ldots
v_1)\,(v_i\ldots v_2)\,(v_{i+1} \ldots v_3)\, ( v_{i+2} \ldots v_4) \cdot$

\vspace{.1in}
$ \ \ \ \underline{\sigma_1 (v_2 v_3 v_1 v_2 \sigma_1  v_2 v_1 v_3 v_2)} \, (v_4 v_3 v_2 v_1)
\ldots  (v_{i+2}  v_{i+1} v_i v_{i-1})\, ( v_{i+3} \ldots v_j)\, (v_{i+2} \ldots v_{j-1})$

\vspace{.1in}
$ = (v_{j-1}\ldots v_{i+2})\,(v_j\ldots v_{i+3})\, (v_{i-1}\ldots v_1)\,(v_i\ldots v_2)\,(v_{i+1}
\ldots v_3)\, ( v_{i+2} \ldots v_4) \cdot$

\vspace{.1in}
$ \ \ \  (\underline{v_2 }v_3 \underline{v_1} v_2 \sigma_1  v_2 v_1 v_3 v_2) \sigma_1 \,
\underline{(v_4 v_3 v_2 v_1) \ldots (v_{i+2} v_{i+1} v_i v_{i-1})}\, ( v_{i+3} \ldots v_j)\,
(v_{i+2} \ldots v_{j-1})$

\vspace{.1in}
$ \stackrel{(\ddagger)}{=}  (v_{j-1}\ldots v_{i+2})\,(v_j\ldots v_{i+3})\, (v_{i-1}\ldots
v_1)\,(v_i\ldots v_2)\,(v_{i+1}
\ldots v_3 v_2 v_1) \, ( v_{i+2} \ldots v_4 v_3 v_2) \cdot$

\vspace{.1in}
$ \ \ \  \sigma_1  \underline{v_2} v_1 \underline{v_3} v_2 \sigma_1 \, (\underline{v_4 \ldots
v_{i+2}})\, (v_3 \ldots v_{i+1})\, (v_2 \ldots v_i) \, (v_1 \ldots v_{i-1})
\, (\underline{ v_{i+3} \ldots v_j}) \, (v_{i+2} \ldots v_{j-1})$

\vspace{.1in}
$ = (v_{j-1}\ldots v_{i+2})\,(v_j\ldots v_{i+3})\, (\underline{ v_{i-1}\ldots v_1)\,(v_i\ldots
v_2)\,(v_{i+1} \ldots v_3 v_2 v_1) \, ( v_{i+2} \ldots v_4 v_3 v_2})  \cdot$

\vspace{.1in}
$ \ \ \   \sigma_1\, (v_2  v_3 v_4 \ldots v_{i+2}  v_{i+3} \ldots v_j)\, (v_1 v_2 \ldots
v_{i+1} v_{i+2} \ldots v_{j-1})\, 
\sigma_1 \,  (v_2 \ldots v_i) \, (v_1 \ldots v_{i-1})$

\vspace{.1in}
$ \stackrel{(\ddagger)}{=}  (v_{j-1}\ldots v_{i+2})\,(v_j\ldots v_{i+3})\, (v_{i-1} v_i v_{i+1}  
v_{i+2})
\ldots (v_1 \underline{v_2 v_3 } v_4)\,(\underline{v_2} v_1 v_3 v_2) \cdot$

\vspace{.1in}
$ \ \ \  \sigma_1   \, (v_2  \ldots v_j)\, (v_1 \ldots v_{j-1})\, 
\sigma_1 \,  (v_2 \ldots v_i) \, (v_1 \ldots v_{i-1}) = \cdots $

\vspace{.1in}
$ = (v_{j-1}\ldots v_{i+2})\,v_{i+1}\,(v_j\ldots v_{i+3})\, (v_{i-1} v_i v_{i+1}   v_{i+2})
\ldots (\underline{v_1 v_2 }v_3 v_4)\,( \underline{v_1} v_3 v_2) \cdot$

\vspace{.1in}
$ \ \ \   \sigma_1  \, (v_2  \ldots v_j)\, (v_1 \ldots v_{j-1})\, 
\sigma_1 \,  (v_2 \ldots v_i) \, (v_1 \ldots v_{i-1}) = \cdots $

\vspace{.1in}
$ = (v_{j-1}\ldots v_{i+2})\,v_{i+1}v_i\,(v_j\ldots v_{i+3})\, (v_{i-1} v_i v_{i+1}   v_{i+2})
\ldots (v_1 v_2 \underline{v_3 v_4})\,(\underline{ v_3} v_2) \cdot$

\vspace{.1in}
$ \ \ \  \sigma_1  \, (v_2  \ldots v_j)\, (v_1 \ldots v_{j-1})\, 
\sigma_1 \,  (v_2 \ldots v_i) \, (v_1 \ldots v_{i-1}) = \cdots $

\vspace{.1in}
$ = (v_{j-1}\ldots v_i)\,(v_j\ldots v_{i+3})\,v_{i+2}\, (v_{i-1} v_i v_{i+1}  
v_{i+2})
\ldots (v_1 \underline{v_2 v_3} v_4)\,\underline{v_2} \cdot$

\vspace{.1in}
$ \ \ \   \sigma_1 \, (v_2  \ldots v_j)\, (v_1 \ldots v_{j-1})\, 
\sigma_1 \,  (v_2 \ldots v_i) \, (v_1 \ldots v_{i-1}) = \cdots $

\vspace{.1in}
$ = (v_{j-1}\ldots v_i)\,(v_j\ldots v_{i+3})\,v_{i+2} v_{i+1} \, \underline{ ( v_{i-1} v_i v_{i+1}  
v_{i+2}) \ldots (v_1 v_2 v_3 v_4)} \cdot$

\vspace{.1in}
$ \ \ \ \sigma_1  \,  (v_2  \ldots v_j)\, (v_1 \ldots v_{j-1})\, 
\sigma_1 \,  (v_2 \ldots v_i) \, (v_1 \ldots v_{i-1}) $

\vspace{.1in}
$\stackrel{(\ddagger)}{=}  (v_{j-1}\ldots v_i)\, (v_j\ldots  v_{i+1})\, \underline{(v_{i-1} 
\ldots v_1)} \, ( v_i \ldots v_2) \,  ( v_{i+1} \ldots v_3)  \,  ( v_{i+2} \ldots v_4)  \cdot$

\vspace{.1in}
$ \ \ \  \underline{ \sigma_1} \, (v_2  \ldots v_j)\, (v_1 \ldots v_{j-1})\, 
\sigma_1 \,  (v_2 \ldots v_i) \, (v_1 \ldots v_{i-1}) $

\vspace{.1in}
$ = (v_{j-1}\ldots v_1)\,(v_j\ldots  v_2)\,\sigma_1\,
 \underline{(v_{i+1} \ldots  v_3) \, ( v_{i+2} \ldots v_4) \, (v_2  \ldots v_j)\, (v_1 \ldots
v_{j-1})}  \cdot$

\vspace{.1in}
$ \ \ \ \sigma_1 \, (v_2\ldots v_i) \, (v_1 \ldots v_{i-1}).$

\bigbreak

\noindent On the other hand we have: 

\vspace{.1in}
\noindent $ \sigma_j \sigma_i \stackrel{(*)}{=}  (v_{j-1}\ldots v_1)\,(v_j\ldots v_2)\,\sigma_1\,
\underline{(v_2 \ldots v_j)\,(v_1 \ldots v_{j-1}) \, (v_{i-1}\ldots v_1)\,(v_i\ldots v_2)}
\cdot $

\vspace{.1in}
 $ \ \ \ \sigma_1 \, (v_2 \ldots v_i)\,(v_1 \ldots v_{i-1}).$

\bigbreak

\noindent Therefore, in order that $ \sigma_i \sigma_j = \sigma_j \sigma_i$ it suffices to show that
the underlined expressions above are equal. Indeed,  we have: 

\vspace{.1in}
\noindent $  (v_{i+1} \ldots  v_3) \, ( v_{i+2} \ldots v_4) \, (\underline{v_2} v_3 \ldots v_j)\,
(v_1 \ldots v_{j-1})$

\vspace{.1in}
$= (v_{i+1} \ldots  v_3 v_2) \, (\underline{ v_{i+2} \ldots v_4) \, (v_3 v_4 \ldots  v_{i+2}} 
v_{i+3} \ldots    v_j)\, (v_1 \ldots v_{j-1})$

\vspace{.1in}
$\stackrel{(\dagger)}{=} ( \underline{v_{i+1} \ldots  v_3 v_2) \, (v_3 \ldots v_{i+1}} v_{i+2}
v_{i+1} \ldots  v_3)\,  (v_{i+3} \ldots    v_j)\, (v_1 \ldots v_{j-1})$

\vspace{.1in}
$\stackrel{(\dagger)}{=} (\underline{ v_2 \ldots v_i v_{i+1}} v_i \ldots  v_2) \, (
\underline{v_{i+2}} v_{i+1}
\ldots  v_3)\,  ( \underline{ v_{i+3} \ldots    v_j})\, (\underline{v_1 } v_2 \ldots v_{j-1})$

\vspace{.1in}
$= (v_2 \ldots  v_{i+1} v_{i+2} v_{i+3} \ldots v_j)\, (v_i \ldots  v_2) \, v_1 \, 
( \underline{ v_{i+1} \ldots v_3)\,  (v_2 v_3 \ldots
v_{i+1}} v_{i+2} \ldots v_{j-1})$

\vspace{.1in}
$\stackrel{(\dagger)}{=} (v_2 \ldots v_j)\, ( \underline{v_i \ldots v_2 v_1) \, 
(v_2 \ldots v_i} v_{i+1} v_i \ldots v_2)\, (v_{i+2} \ldots v_{j-1})$

\vspace{.1in}
$\stackrel{(\dagger)}{=} (v_2 \ldots v_j)\, (\underline{ v_1 \ldots v_{i-1} v_i } v_{i-1} \ldots
v_1) \,  (\underline{v_{i+1}} v_i  \ldots v_2)\, ( \underline{v_{i+2} \ldots v_{j-1} })$

\vspace{.1in}
$= (v_2 \ldots v_j)\, (v_1 \ldots  v_{j-1}) \,  ( v_{i-1} \ldots v_1)\, ( v_i  \ldots v_2).$ 
   \hfill $\Box$

\bigbreak

 By Lemmas 1, 2, and 3  the proof of Theorem 2 is now concluded.  QED.

\section{A Reduced Presentation for the Flat Virtual Braid Group}

The flat virtual braids were introduced in \cite{SVKT}. As with the virtual braids, the set of flat
virtual braids on $n$ strands forms a group, the {\it flat virtual braid group,}  denoted $FV_n.$ 
The generators of $FV_n$ are the virtual crossings $v_{1},\ldots, v_{n-1}$ and the flat crossings
$c_{1},\ldots, c_{n-1},$ which --as already said in 1.1-- can be seen as immersed crossings. See Figure
19. 

\smallbreak 
 So, flat crossings and  virtual crossings both represent geometrically the generators of the
symmetric group
$S_n.$  But the mixed relations between them are not symmetric (see below). In fact, the flat virtual
braid group is the quotient of the virtual braid group $VB_{n}$ modulo the relations ${\sigma_i}^2  = 1$
for all $i.$ Thus,  $FV_n$ is the free product of two copies of $S_n,$ modulo the set of mixed relations
specified below. Note that $FV_1 = S_2 * S_2$ (no extra relations), and it is infinite.
\bigbreak

Recall that in section $1.1$ we discussed flat virtual knots and links, and that we pointed out that this category is equivalent to
the category of virtual strings developed in \cite{TURAEV}. Just so, the flat virtual braids are the appropriate theory of braids for
the category of virtual strings. Every virtual string is the closure of a flat virtual braid.

\bigbreak

 $$\vbox{\picill10cmby1.5in(VP19)  }$$

\begin{center}
{\bf Figure 19 -- The Generators of $FV_n$}
\end{center}
\vspace{3mm}

The virtual generators satisfy among themselves the  virtual
relations.  Similarly, the flat  generators satisfy among themselves the following {\it
flat relations:}  

$$ \begin{array}{cccl} 
{c_{i}}^2 & = & 1, &   \\
c_{i}c_{i+1}c_{i} & = & c_{i+1}c_{i} c_{i+1}, &   \\
c_i c_j & = & c_j c_i, & \mbox{for} \ j\neq i\pm 1.  \\
\end{array}$$

\noindent The {\it mixed flat relations} between flat and virtual generators are as follows:
 
$$ \begin{array}{cccl} 
c_i v_j & = & v_j c_i, & \mbox{for} \ j\neq i\pm 1,  \\
v_{i} c_{i+1} v_{i} & = & v_{i+1} c_i v_{i+1}. &   \\ 
\end{array}$$

 The second mixed relation will be called the {\it special detour flat relation} and it is
illustrated in Figure 9. 
Then, as for the virtual braids, we have for the flat crossings the inductive defining relations:
$c_{i+1} = v_i v_{i+1} c_i v_{i+1} v_i,$ which leads to the {\it defining relations:} 

$$c_{i+1} := (v_i\ldots v_2v_1)\,(v_{i+1}\ldots v_3v_2)\,c_1\, (v_2v_3\ldots
v_{i+1})\,(v_1v_2\ldots v_i) $$

\bigbreak

\noindent for \  $i=2,\ldots, n-1.$ In terms of flat braid diagrams, this relation is the flat braid
detour move of the strands $1, 2, \ldots, i-1$  around the flat crossing $c_i.$ In complete analogy to
the virtual braid group we now have  the following:

\begin{th}{ \ The flat virtual braid group $FV_{n}$ has the following reduced presentation.  

\[FV_{n} = \left< \begin{array}{ll}  
\begin{array}{l} 
c_1, v_1, \ldots ,v_{n-1}  \\
\end{array} &
\left|
\begin{array}{l} 

v_i v_{i+1} v_i = v_{i+1} v_i v_{i+1}    \\

v_i v_j=v_j v_i, \ \  \ j\neq i\pm 1   \\ 

{c_1}^2  = 1,  \ {v_i}^2  = 1, \ \ \ 1 \leq i \leq n-1  \\

c_1 v_j=v_j c_1, \ \ \ j>2   \\ 

 (v_1 c_1 v_1)\,(v_2 c_1 v_2)\, (v_1 c_1 v_1) =
  (v_2 c_1 v_2)\, (v_1 c_1 v_1) \, (v_2 c_1 v_2)   \\ 

c_1\, (v_2 v_3 v_1 v_2 c_1 v_2 v_1 v_3 v_2) =
    (v_2 v_3 v_1 v_2 c_1 v_2 v_1 v_3 v_2) \, c_1   \\ 
\end{array} \right.  \end{array} \right>.  \]
}
\end{th}

\section{Other Categories}

Welded  braids were introduced in \cite{FRR}. They satisfy the same isotopy relations as 
the virtuals, but for welded braids one also allows one of the two forbidden moves, the move $(F_1)$ of
Figure~4, which  contains an over arc and one virtual crossing. One can consider welded knots and
links in this way, and the explanation for the choice of moves lies in the fact that the first forbidden
move preserves the  combinatorial fundamental group. This not true
for the other forbidden move
$(F_2).$ The corresponding  {\it welded braid group}   on $n$ strands, 
$WB_n,$ has the same generators and relations as the virtual braid group, but with the extra relations:

$$v_{i} \sigma_{i+1} \sigma_{i}  =  \sigma_{i+1} \sigma_i v_{i+1}   \  \  \ (F_1)$$

\noindent Figure 4 illustrates a variant of this relation. Just as in the virtual braid group, the
braiding generators $\sigma_2, \ldots, \sigma_{n-1}$ of the welded braid group can be written in terms of
$\sigma_1$ and the welded generators $v_1,\ldots,v_{n-1}.$  We can then obtain a reduced presentation 
for  $WB_n$  with generators 

$$\{\sigma_1, v_1,\ldots,v_{n-1}\}$$ 

\noindent and the defining relations:

$$\sigma_{i+1} := (v_i\ldots v_2v_1)\,(v_{i+1}\ldots v_3v_2)\,\sigma_1\, (v_2v_3\ldots
v_{i+1})\,(v_1v_2\ldots v_i)  \  \  \ (*)$$

\noindent By the box detour trick (see Remark 2) we can easily  reduce the set of extra
relations $(F_1)$ to the basic relation:

$$v_{1} \sigma_{2} \sigma_{1}  =  \sigma_{2} \sigma_1 v_{2}, $$

\noindent  which with the substitution $ \sigma_{2}  =   v_{1} v_{2} \sigma_1 v_{2}  v_{1} $ is
equivalent to the relation: 

$$ v_{2} \sigma_1 v_{2}  v_{1} \sigma_{1}  =  v_{1} v_{2} \sigma_1 v_{2}  v_{1} \sigma_1
v_{2}. $$

\noindent  Thus we have:

\begin{th}{ \ The welded braid group $WB_{n}$ has the following reduced presentation.  

\[WB_{n} = \left< \begin{array}{ll}  
\begin{array}{l} 
\sigma_1, v_1, \ldots ,v_{n-1}  \\
\end{array} &
\left|
\begin{array}{l} 

v_i v_{i+1} v_i = v_{i+1} v_i v_{i+1}    \\

v_i v_j=v_j v_i, \ \  \ j\neq i\pm 1   \\ 

{v_i}^2  = 1, \ \ \ 1 \leq i \leq n-1  \\

\sigma_1 v_j=v_j \sigma_1, \ \ \ j>2   \\ 

 (v_1 \sigma_1 v_1)\,(v_2 \sigma_1 v_2)\, (v_1 \sigma_1 v_1) =
  (v_2 \sigma_1 v_2)\, (v_1 \sigma_1 v_1) \, (v_2 \sigma_1 v_2)   \\ 

 v_{1} \, (v_{2} \sigma_1 v_{2}  v_{1} \sigma_{1})  =   (v_{2} \sigma_1 v_{2}  v_{1} \sigma_1)\, 
v_{2}    \\ 

\sigma_1\, (v_2 v_3 v_1 v_2 \sigma_1 v_2 v_1 v_3 v_2) =
    (v_2 v_3 v_1 v_2 \sigma_1 v_2 v_1 v_3 v_2) \, \sigma_1   \\ 
\end{array} \right.  \end{array} \right>.  \]
}
\end{th}

Note now that the relations $(F_1)$ can be regarded as a way of detouring sequences of classical 
crossings over welded crossings, via the  inductive relations: 

$$ v_{i+1}   = {\sigma_{i}}^{-1} {\sigma_{i+1}}^{-1} v_i \sigma_{i+1} \sigma_i,$$

\noindent which, by induction, lead to the defining relations:

$$v_{i+1} := ({\sigma_i}^{-1}  \ldots {\sigma_i}^{-1})\,({\sigma_{i+1}}^{-1} \ldots
{\sigma_2}^{-1})\,v_1\, (\sigma_2 \ldots \sigma_{i+1})\,(\sigma_1 \ldots \sigma_i)  \  \  \ (**)$$

\bigbreak

\noindent for \  $i=1,\ldots, n-2.$ By the box detour trick we reduce the relations involving welded
generators. For example, the welded relations reduce to the following two basic ones:
$$v_1 v_2 v_1 = v_2 v_1 v_2 \ \ \mbox{and} \ \ v_1 v_3 = v_3 v_1.$$

\noindent Thus, we obtain the following reduced presentation for $WB_{n} $ with a single welded
generator.

\begin{th}{ \ The welded braid group $WB_{n}$ has the following reduced presentation.  

\[WB_{n} = \left< \begin{array}{ll}  
\begin{array}{l} 
v_1, \sigma_1, \ldots ,\sigma_{n-1}  \\
\end{array} &
\left|
\begin{array}{l} 

\sigma_i \sigma_{i+1} \sigma_i = \sigma_{i+1} \sigma_i \sigma_{i+1}    \\

\sigma_i \sigma_j=\sigma_j\sigma_i, \ \  \ j\neq i\pm 1   \\ 

{v_1}^2  = 1  \\

v_1 \sigma_j=\sigma_j v_1, \ \ \ j>2   \\ 

 (\sigma_1 v_1 {\sigma_1}^{-1})\, ( {\sigma_2}^{-1} v_1 \sigma_2) \, (\sigma_1 v_1 {\sigma_1}^{-1}) = \\ 

  ( {\sigma_2}^{-1} v_1 \sigma_2) \, (\sigma_1 v_1 {\sigma_1}^{-1}) \, ( {\sigma_2}^{-1} v_1 \sigma_2) 
\\ 

v_1\, ({\sigma_2}^{-1} {\sigma_1}^{-1} {\sigma_3}^{-1} {\sigma_2}^{-1} v_1 \sigma_2 \sigma_3 \sigma_1
\sigma_2) = \\

({\sigma_2}^{-1} {\sigma_1}^{-1} {\sigma_3}^{-1} {\sigma_2}^{-1} v_1 \sigma_2 \sigma_3
\sigma_1 \sigma_2) \, v_1   \\ 
\end{array} \right.  \end{array} \right>.  \]
}
\end{th}

Another generalization of the virtual braid group is obtained by adding both types of forbidden moves
(recall Figure 4). We call this the {\it unrestricted virtual braid group,} denoted $UB_n.$ It is known
that any classical knot can be unknotted in the virtual category if we allow both forbidden moves
\cite{KANENOBU, NELSON}. Nevertheless, linking phenomena still remain.  The unrestricted braid group
itself is non trivial, deserving further study. 

\bigbreak
\noindent By adding both types of forbidden moves: 

$$v_{i} \sigma_{i+1} \sigma_{i}  =  \sigma_{i+1} \sigma_i v_{i+1}  \  \ (F_1) \ \ \ \mbox{and} \ \
\  \sigma_i \sigma_{i+1} v_{i}  =  v_{i+1}  \sigma_{i} \sigma_{i+1}  \  \ (F_2) $$

\noindent and using the defining relations $ (*)$ we obtain a reduced presentation 
for  $UB_n$  with generators $\sigma_1, v_1,\ldots,v_{n-1},$  which is in fact a quotient of the 
corresponding reduced presentation of $WB_n$ by the second forbidden move:

$$ \sigma_1 \sigma_2 v_1  =  v_2  \sigma_1 \sigma_2. $$

\begin{th}{ \ The unrestricted virtual braid group  has the following reduced presentation.  

\[UB_{n} = \left< \begin{array}{ll}  
\begin{array}{l} 
\sigma_1, v_1, \ldots ,v_{n-1}  \\
\end{array} &
\left|
\begin{array}{l} 

v_i v_{i+1} v_i = v_{i+1} v_i v_{i+1}    \\

v_i v_j=v_j v_i, \ \  \ j\neq i\pm 1   \\ 

{v_i}^2  = 1, \ \ \ 1 \leq i \leq n-1  \\

\sigma_1 v_j=v_j \sigma_1, \ \ \ j>2   \\ 

 (v_1 \sigma_1 v_1)\,(v_2 \sigma_1 v_2)\, (v_1 \sigma_1 v_1) =
  (v_2 \sigma_1 v_2)\, (v_1 \sigma_1 v_1) \, (v_2 \sigma_1 v_2)   \\ 

 v_{1} \, (v_{2} \sigma_1 v_{2}  v_{1} \sigma_{1})  =   (v_{2} \sigma_1 v_{2}  v_{1} \sigma_1)\, 
v_{2}    \\

  (\sigma_1  v_{1} v_{2} \sigma_1 v_{2}) \, v_{1}  =  v_2 \, (\sigma_1  v_{1} v_{2} \sigma_1 v_{2}) 
\\

\sigma_1\, (v_2 v_3 v_1 v_2 \sigma_1 v_2 v_1 v_3 v_2) =
    (v_2 v_3 v_1 v_2 \sigma_1 v_2 v_1 v_3 v_2) \, \sigma_1   \\

\end{array} \right.  \end{array} \right>.  \]
}
\end{th}

Just as in the case of welded braids, we can also give a reduced presentation with one virtual generator
and $n-1$ braiding generators. For unrestricted virtual braids, there are two possible such reduced
presentations, depending upon using either the first or the second forbidden move in performing the
detour substitutions. In the $(F_1)$ case the defining relations are given by (**) and, thus, the reduced
presentation is a quotient of the corresponding presentation of the welded braid group by the following
relation:

$$\sigma_1 \sigma_2 v_{1} {\sigma_2}^{-1} {\sigma_1}^{-1} =   {\sigma_1}^{-1}{\sigma_2}^{-1}v_{1} 
\sigma_2  \sigma_1. $$

\noindent Since in this presentation $v_2$ is defined via $(F_1)$ in terms of $v_1,$ the reader will
 note that the transcription of this last relation appears to be a mixture of $(F_1)$ and $(F_2).$
Similarly we could have started with $(F_2)$ and obtained first an analogue of the welded braid group
and then, adding $(F_1),$ obtained another reduced presentation of the unrestricted virtual braid group.

\bigbreak

Finally, we define the {\it flat unrestricted braid group,} denoted $FU_{n},$ to be the quotient
of the flat virtual braid group $FV_{n}$ (see Theorem 3) by the forbidden moves of $FV_{n}$ (see
Figure 9):

$$c_{i} c_{i+1} v_{i}  =  v_{i+1} c_i c_{i+1}.$$

\noindent Note that for the flat virtual braid group there is only one type of forbidden move.

\begin{th}{ \ The flat unrestricted  braid group $FU_{n}$ has the following reduced
presentation.  

\[FU_{n} = \left< \begin{array}{ll}  
\begin{array}{l} 
c_1, v_1, \ldots ,v_{n-1}  \\
\end{array} &
\left|
\begin{array}{l} 

v_i v_{i+1} v_i = v_{i+1} v_i v_{i+1}    \\

v_i v_j=v_j v_i, \ \  \ j\neq i\pm 1   \\ 

{c_1}^2  = 1, \ {v_i}^2  = 1,  \ \ \ 1 \leq i \leq n-1  \\

c_1 v_j=v_j c_1, \ \ \ j>2   \\ 

 (v_1 c_1 v_1)\,(v_2 c_1 v_2)\, (v_1 c_1 v_1) =
  (v_2 c_1 v_2)\, (v_1 c_1 v_1) \, (v_2 c_1 v_2)   \\ 

 v_{1} \, (v_{2} c_1 v_{2}  v_{1} c_{1})  =   (v_{2} c_1 v_{2}  v_{1} c_1)\, 
v_{2}    \\

c_1\, (v_2 v_3 v_1 v_2 c_1 v_2 v_1 v_3 v_2) =
    (v_2 v_3 v_1 v_2 c_1 v_2 v_1 v_3 v_2) \, c_1   \\ 
\end{array} \right.  \end{array} \right>.  \]
}
\end{th}

\begin{rem}{ \rm  \  Note that the flat unrestricted braid group $FU_{n}$ is a free product with
amalgamation of two copies of the symmetric group $S_n.$ An unreduced presentation of  $FU_{n}$ can de
configured to be symmetric with respect to the roles of the generators $c_i$ and $v_i.$ As a result
there is another reduced presentation that can be obtained from the reduced presentation above by
interchanging the roles of  $v_i$ and $c_i.$ 
 }
\end{rem}

\begin{rem}{ \rm  \  Note that the flat unrestricted braid group $FU_{n}$ is also a quotient of the
welded braid group $WB_{n}$ (see Theorem 4),  obtained by setting all the squares of the braiding
generators equal to 1. Thus there is a surjective homomorphism from $WB_{n}$  to $FU_{n}.$ This
homomorphism is a direct analogue of the standard homomorphism from $B_n$ to the symmetric group $S_n.$ 
Figure 20 gives a commutative diagram of these relationships. Note that all structures map eventually to the 
symmetric group $S_{n}.$ In the case of the virtual braids and their quotients, this map to the symmetric group
takes the same value on virtual generators $v_{i}$ and braiding generators $\sigma_{i}.$ The itermediate
mappings to $FV_{n}$ and $FU_{n}$ preserve these differences.}
\end{rem}

\vspace{3mm}

{\tt    \setlength{\unitlength}{0.92pt}
\begin{picture}(387,139)
\thicklines   \put(228,55){\vector(-3,-1){110}}
              \put(149,63){\vector(1,0){69}}
              \put(257,130){\vector(1,0){69}}
              \put(145,131){\vector(1,0){69}}
              \put(33,133){\vector(1,0){69}}
              \put(15,124){\vector(3,-4){83}}
              \put(339,123){\vector(-2,-1){96}}
              \put(111,55){\vector(0,-1){39}}
              \put(230,116){\vector(0,-1){39}}
              \put(112,116){\vector(0,-1){39}}
\thinlines    \put(101,4){$S_{n}$}
              \put(222,59){$FU_{n}$}
              \put(104,60){$FV_{n}$}
              \put(330,127){$UB_{n}$}
              \put(217,127){$WB_{n}$}
              \put(103,130){$VB_{n}$}
              \put(1,130){$B_{n}$}
\end{picture}}

\begin{center}
{\bf Figure 20 -- A Diagram of Relationships}
\end{center}

\section{Welded Braids and Tubes in Four-Space}

The welded braid group $WB_{n}$ can be interpreted as the fundamental group of the configuration space of
 $n$ disjoint circles trivially embedded in three dimensional space $\RR^{3}$. This group (the so-called motion
group of disjoint circles) can, in turn, be interpreted as a braid group of tubes imbedded in $\RR^{3} \times \RR
= \RR^{4}.$ These braided tubes in four-space are generated by two types of elementary braiding. In Figure 21, we show diagrams that
can be interpreted as immersions of tubes in three-space. Each such immersion is a projection of a corresponding embedding in
four-space. The first two diagrams of Figure 21 each illustrate a tube passing through another tube. When tube $A$ passes through
tube $B$ we make a corresponding classical braiding crossing with arc $A$ passing under arc $B.$ The four-dimensional
interpretation of tube $A$ passing through tube $B$ is that:
 As one looks at the levels of intersection with 
$\RR^{3} \times t$ for different values of $t,$ one sees two circles $A(t)$ and $B(t).$ As the variable $t$ 
increases, the $A(t)$ circle (always disjointly embedded from the $B(t)$ circle) moves through the $B(t)$
circle. This process is illustrated in Figure 22.
\bigbreak

While the classical crossing in a welded braid diagram
corresponds to a genuine braiding of the tubes in four-space (as described above), the virtual crossing 
corresponds to tubes that do not interact in the immersion representation  (see again Figure 21). These non-interacting tubes can pass
over or under each other, as these local projections correspond to equivalent embeddings in four-space.
\bigbreak

$$\vbox{\picill3inby5in(VP21)  }$$

\begin{center}
{\bf Figure 21 -- Tubular Correspondences}
\end{center}
\vspace{3mm}

$$\vbox{\picill4inby2.6in(VP22)  }$$

\begin{center}
{\bf Figure 22 -- Braiding of Circles}
\end{center}
\vspace{3mm}

It is an interesting exercise to verify that the moves in the welded braid group each induce equivalences of the
corresponding tubular braids in four-space. In particular, the move $(F_1)$ induces such an isotopy, while the forbidden move
$(F_2)$ does not. For more on this subject, the reader can consult
\cite{DVK} and the references therein.
\bigbreak

Consider now the surjection $WB_{n} \longrightarrow FU_{n}$ from the welded braid group to the flat unrestricted
braids. Flat unrestricted braids can be represented by immersions of tubes in three-space as illustrated also in
Figure 21. There we have indicated a decorated immersion of two intersecting tubes as the correspondent of the
flat classical crossing in $FU_{n}.$ One must specify the rules for handling these immersions in order
to obtain the correspondence. We omit that discussion here, but point out the interest in having a uniform
context for the surjection of the welded braids to the flat unrestriced braids. The flat unrestricted braids
carry the distinction between braided flat and welded flat crossings and otherwise keep track of the relative 
permutations of these two types of crossing.
\bigbreak

\bigbreak

\bigbreak

\noindent {\bf Acknowledgement.} Much of this effort was sponsored by the Defense
Advanced Research Projects Agency (DARPA) and Air Force Research Laboratory, Air
Force Materiel Command, USAF, under agreement F30602-01-2-05022. 
The U.S. Government is authorized to reproduce and distribute reprints
for Government purposes notwithstanding any copyright annotations thereon. The
views and conclusions contained herein are those of the author and should not be
interpreted as necessarily representing the official policies or endorsements,
either expressed or implied, of the Defense Advanced Research Projects Agency,
the Air Force Research Laboratory, or the U.S. Government. (Copyright 2004.) 
It gives the first author great pleasure to acknowledge support from NSF Grant DMS-0245588,
and to give thanks to the University of Waterloo and the Perimeter Institute in Waterloo, 
Canada for their hospitality during the preparation of this research.

%%%%%%%%%%%%%%%%%%%%%%%%%%%%%%%%%%%%%%%%%%%%%%%%%%%%%%%%%%%%%%%%%%%%%%%%%%%%%%%
\bigbreak

\noindent {\sc L.H. Kauffman: Department of Mathematics, Statistics and
Computer Science, University of Illinois at Chicago, 851 South Morgan St., Chicago IL 60607-7045, USA. }

 \vspace{.1in}
 
\noindent {\sc S. Lambropoulou: Laboratoire de Math\'ematiques Nicolas Oresme,  Universit\'e de
Caen, F-14032 Caen cedex, France and 

\noindent Department of Mathematics, 
National Technical University of Athens,
Zografou campus, GR-157 80 Athens, Greece.}

\vspace{.1in}
\noindent {\sc E-mails:} \ {\tt 
kauffman@uic.edu  \ \ \  sofia@math.ntua.gr and sofia@math.unicaen.fr}

\noindent {\sc URLs:} \ {\tt 
www.math.uic.edu/~kauffman  \ \ \ \ \ \ \ \ http://users.ntua.gr/sofial}

%%%%%%%%%%%%%%%%%%%%%%%%%%%%%%%%%%%%%%%%%%%%%%%%%%%%%%%%%%%%%%%%%%%%%%%%%%%%%%%

\end{document}